\documentstyle[12pt]{article}

\newtheorem{definition}{Definition}[section]
\newtheorem{theorem}[definition]{Theorem}
\newtheorem{proposition}[definition]{Proposition}

\newtheorem{lemma}[definition]{Lemma}

\newcommand{\qed}{\hfill \rule{.5em}{1em}}

\begin{document}

\title{The Euler Series of Restricted Chow Varieties}
\author{Javier Elizondo}
\date{{\small  8 April 1993}}

\maketitle

We are interested in understanding the
geometry and topology of Chow varieties and the space of
all effective cycles  with the same homology class,
``restricted Chow variety''.
In this paper we discuss the case of
projective algebraic varieties on which an algebraic torus acts linearly
such that there are finitely many irreducible invariant subvarieties.
The projective toric varieties are important examples which
satisfy this condition.

Let $X$
be a  projective  algebraic variety together with a fixed
embedding.
Denote by ${\cal C}_{\lambda}$
the space of all effective cycles on $X$ whose homology class
is $\lambda.$ This space is in fact a subvariety of a Chow variety of $X.$
Observe that if $X \, = \, {\bf P}^{n}$ and
$\lambda \, = \, d[{\bf P}^{p}]$ we have the following equality
$$
{\cal C}_{\lambda}({\bf P}^{n}) \, = \,
{\cal C}_{p,d}({\bf P}^{n})
$$
where ${\cal C}_{p,d}({\bf P}^{n})$ is the Chow variety of $X$ which
parametrizes the effective cycles of ${\bf P}^{n}$ of dimension
$p$ and degree $d.$
 Very little is known about ${\cal C}_{p,d}(X)$.
In  \cite{la&yau-hosy} H.B.Lawson and S.S.Yau posed the problem of
computing the Euler characteristic of
${\cal C}_{\lambda}.$ They introduced a formal power series
which satisfies, for $X \, = \, {\bf P}^{n},$ the following equality
\begin{equation}
\label{layau}
\sum_{d=o}^{\infty} \chi \left({\cal C}_{p,d} ({\bf P}^{n})\right) t^{d}
\, = \, {\left( \frac{1}{1-t} \right)}^{(_{p+1}^{n+1})}.
\end{equation}
This shows that the formal power series on the left is rational and
solves the problem of computing the Euler characteristic of
${\cal C}_{p,d} ({\bf P}^{n}).$ They also computed the case
$X\, = \, {\bf P}^{n} \times {\bf P}^{m}.$
When we try to generalize the above series in order to compute
$\chi ( {\cal C}_{\lambda})$ we face some problems. If a
basis for the integral homology group (modulo torsion) of $X$ is fixed,
and we write it as a multiplicative group in the canonical way,
we arrive at a series where some of the powers could be negative integers.
In others words, we do not obtain a formal power series, and
the concept of rationality is not defined a priori. In the author's thesis
\cite{eli-thesis} an ``ad hoc'' definition of
rationality was given and it was proved
that  for any smooth projective toric variety
the corresponding series is rational. In this paper we follow an
intrinsic formulation
suggested by E.Bifet \cite{emili-cartas} in order to reinterpret
the results in \cite{eli-thesis} and  prove them for
any projective algebraic variety with a finite number of invariant
subvarieties
under an algebraic torus action. The definition of the monoid $C,$
rationality and lemma \ref{fibers} are taken from him.
The main results can be stated in
the following way.
Let $C$ be the monoid of homology classes of effective p-cycles, and
let ${\bf Z}[[C]]$ be the ring of functions (with respect to
the convolution product) over $C.$ Denote by  ${\bf Z}[C]$ the ring
of functions with finite support on $C.$ We say that an element of
${\bf Z}[[C]]$ is rational if it is the quotient of two elements
in ${\bf Z}[C].$
We define the Euler series of $X$ by
\begin{equation}
\label{euler}
E_{p} \, = \, \sum_{\lambda\in{C}} \, \chi \left({\cal C}_{\lambda}\right)
\lambda
  \, \, \, \in \, {\bf Z}[[C]]
\end{equation}
This generalizes the series in equation (\ref{layau}).
 Let $X$ be a projective algebraic variety and $V_1 , \ldots , V_N$ be
its p-dimensional invariant subvarieties under the action of
the algebraic torus. Denote by
 $e_{[V_{i}]} \, \in \, {\bf Z}[C]$
the characteristic function of the set
$\{[V_{i}]\}.$
 Theorem
{}~\ref{main} says,
\begin{equation}
\label{int}
E_{p} \, = \,
\prod_{i=1}^{N} \left( \frac{1}{1-e_{[V_{i}]}} \right)
\end{equation}
If $X$ is a smooth projective toric variety we define $C_{T}$
as the monoid of equivariant homology classes of invariant
p-cycles. Denote by  ${\bf Z}[[C_{T}]]$ and by ${\bf Z}[C_{T}]$
the ring of functions  and the ring of functions
with finite support on  $C_{T}$ respectively. An element
of  ${\bf Z}[[C_{T}]]$ is rational if it is the quotient of two elements
of  ${\bf Z}[C_{T}].$
In section three is defined the equivariant Euler
series $E^{T}_{p}$
and it is proved in theorem ~\ref{ec-eec} that
$E^{T}_{p}$ is rational. It is also obtained an explicit formula
for it. Furthermore, a ring homomorphism
$$
J: {\bf Z}[[C_{T}]] \longrightarrow {\bf Z}[[C]]
$$
from  ${\bf Z}[[C_{T}]]$ to  ${\bf Z}[[C]]$ is defined and
we recover formula (\ref{int}) from the following equality
$$
J\left( E_{p}^{T}\right) \, = \, E_{p}.
$$

{\bf Acknowledgement:} This paper is based on the author's Ph.D. thesis
written under the supervision of Blaine Lawson at Stony Brook.
I would like to thank
both Blaine Lawson and Emili Bifet for many wonderful conversations
and ideas they shared with me. I also would like to thank
Paulo Lima-Filho, Santiago Lopez de Medrano and E. Bifet for reading
carefully a first version of this paper and for many
rich comments.

\section{The Euler Series}

In  \cite{la&yau-hosy} B.Lawson and S.S.Yau introduced a
series that becomes a formal power series when a basis for
homology is fixed. They proved it is a rational function for the cases of
 ${\bf P}^{n}$ and   ${\bf P}^{n} \times {\bf P}^{m}.$
In this section we define a more general series
and state
the problem of its rationality in intrinsic terms
 for any algebraic projective variety.
 We follow an approach suggested by E.Bifet  \cite{emili-cartas}.
The definition of the monoid $C,$ rationality
 and lemma ~\ref{fibers} are taken from him.
We start with some basic definitions and some of their properties.

Throughout this article any algebraic projective variety $X$ comes with
a fixed embedding  $ X \stackrel{j}{\hookrightarrow} {\bf P}^{n}.$
An {\bf effective p-cycle} c on $X$ is a finite (formal) sum
\mbox{$c=\sum n_{s} \, V_{s}$} where  each $n_{s}$ is a
non-negative integer and each $V_{s}$  is an irreducible
\mbox{p-dimensional} subvariety of $X$. From now on,
 we shall use the term {\bf cycle} for {\bf effective cycle}.
For any  projective algebraic variety
 $ X \stackrel{j}{\hookrightarrow} {\bf P}^{n}$
we denote by ${\cal C}_{p,d} (X)$ the {\bf Chow variety} of $X$ of
 all cycles of dimension p and degree d in
 ${\bf P}^{n}$ with support on $X.$
By convention, we write ${\cal C}_{p,0} (X) \, = \, \{ \emptyset \}.$
Let $\lambda$ be an element in  $H_{2p} (X,\bf {Z})$ and denote by
  ${\cal C}_{\lambda} (X)$ the space of all
cycles on $X$ whose homology class is  $\lambda.$
Note that ${\cal C}_{\lambda} (X)$ is contained in ${\cal C} _{p,d} (X),$
 where
 $j_{\ast} \, \lambda \, = \, d[{\bf {P}}^{p}].$
\begin{lemma}
\label{compacto}
Let $\lambda$ be an element of $H_{2p}(X,{\bf Z})$, then
 ${\cal C}_{\lambda} (X)$ is a projective algebraic variety.
\end{lemma}
{\bf Proof:}
Let us write
 ${\cal C} _{p,d} (X)= \cup _{i=1}^{M} {\cal C} _{p,d}^{i} (X)$,
  where ${\cal C} _{p,d}^{i} (X)$
are its irreducible components.
Suppose
 ${\cal C}_{\lambda} (X) \cap {\cal C} _{p,d}^{i_{o}} (X) \neq \emptyset.$
 Any two cycles in
${\cal C} _{p,d}^{i_{o}}$ are algebraically equivalent,
 hence they represent the same element in homology. Therefore
 ${\cal C} _{p,d}^{i_{o}} (X) \subset {\cal C}_{\lambda} (X)$
 for some $\lambda.$
Consequently ${\cal C}_{\lambda} (X)=
\cup_{j=1}^{l} {\cal C} _{p,d}^{i_{j}} (X),$
 where
 ${\cal C}_{\lambda} (X) \cap {\cal C} _{p,d}^{i_{j}} (X) \neq \emptyset$
for $j \, = \, 1, \ldots ,l .$
\qed

 For an effective $p$-cycle $c \,=\, \sum \, n_i \, V_i,$ we denote
 by $[c]$ its homology class in $H_{2p}\,(X,{\bf Z}).$
Now, let $C$ be {\bf the monoid of homology classes of effective p-cycles}
in $H_{2p}\,(X,{\bf Z}),$ and let ${\bf Z}[[C]]$ be
 the {\bf set of all integer valued functions on $C$}. We shall write the
 elements of ${\bf Z}[[C]]$ as
 $$
 \sum_{\lambda\in{C}} \;a_{\lambda}\, \lambda \:\:\:\:\:\mbox{where
}\,\,a_{\lambda}
  \, \in \, {\bf Z}.
 $$
 The following lemma allows us
 to prove that
${\bf Z}[[C]]$ is a ring,

\begin{lemma}
\label{fibers}

Let $C$ be the monoid of homology classes of effective p-cycles on $X.$ \,Then
$$
+: \,C \, \times \, C \longrightarrow \, C
$$
has finite fibers.

\end{lemma}

{\bf Proof:}
If $X \, = \, {\bf P}^{n}$
the result is obvious since $C$ is isomorphic to ${\bf N}.$
Let
$$
j_{\ast} :\, H_{2p}\,(X,{\bf Z}) \, \longrightarrow \, H_{2p}\,({\bf
P}^{n},{\bf Z})
$$
be the homomorphism induced by the embedding
$ X \, \stackrel{j}{\hookrightarrow} \, {\bf P}^{n}. $\,
Denote by $C^{\prime} \, \cong \,  {\bf N}$ the monoid of homology classes of
p-cycles on ${\bf P}^{n}.$ It follows from the proof of
lemma ~\ref{compacto} that
$$
j_{\ast}| _{\scriptstyle C} :\, C \, \longrightarrow \, C^{\prime}
$$
has finite fibers. Finally, the lemma follows from the commutative diagram
below
$$
\begin{array}{ccccc}
&& {\scriptstyle +} && \\
 & C\times C & \longrightarrow & C & \\
{j_{\ast}\times j_{\ast} \,|}_{\scriptstyle {C\times C} }& \downarrow &&
\downarrow &

{j_{\ast}\,|}_{\scriptstyle C} \\
& C^{\prime}\times C^{\prime} & \to & C^{\prime} &  \\
&& {\scriptstyle +} &&
\end{array}
$$
\qed
\\
Directly from this lemma we obtain:

\begin{proposition}
 \label{ring}
 Let $C$ and $Z[[C]]$ defined as above. Then $Z[[C]]$ is a ring under the
convolution
product,
i.e.
$$
(f \cdot g) \,(\lambda ) \,= \,
\sum_{\lambda={\mu_{1}}+{\mu_{2}}} \,f(\mu_{1}) g(\mu_{2}).
$$
 \end{proposition}

We are ready for the following definition.
\begin{definition}
\label{ec}
Let $X$ be a projective algebraic variety. The  {\bf Euler  series}
 of $X,$ in dimension p,
is the element
$$
E_{p} \, =  \sum_{\lambda} \chi  \, ({\cal C}_{\lambda}) \,
 \lambda \, \, \, \; \in \, \, \, \;  {\bf {Z}}[[C]] \; \;
$$
where ${\cal C}_{\lambda}(X)$ is the space of all
effective cycles on $X$ with homology class $\lambda,$ and
 $\chi \, ({\cal C}_{\lambda}(X))$
is the
Euler characteristic of  ${\cal C}_{\lambda}$.
\end{definition}
By convention, if ${\cal C}_{\lambda}$ is the empty set then its Euler
characteristic is zero.

Let ${\bf Z}[C]$ be the monoid-ring  of $C$ over ${\bf Z}.$ This
 ring consists of all
elements of ${\bf Z}[[C]]$ with finite support. We arrive at
the following definition,

\begin{definition}
\label{rational}
We say that an element of ${\bf Z}[[C]]$ is rational if it is
the quotient of two elements of ${\bf Z}[C].$
\end{definition}

{\bf Remark:}
Denote by $H$ the homology group $H_{2p}(X,{\bf Z})$
together with a fixed basis $\cal A.$ Consider $H$ as a multiplicative
group in the standard way and suppose that $C$ is isomorphic to the
monoid of the natural numbers. Then it is easy to see that
${\bf Z}[[C]]$ is isomorphic to the ring of formal power series
in as many variables as the rank of $H.$ Therefore ${\bf Z}[C]$ is
the ring of polynomials and the last definition just says
when a formal power series is a rational function.

We are interested in the following problem.\newline
{\bf Problem:}
\label{problem}
\, When is the {\bf Euler series} rational in
the sense of the last definition?

For cycles of dimension zero on a projective algebraic variety, we have
directly from the computation in \cite{mcd-sym} that
$$
E_{0} \, = \, \frac{1}{{\left( 1-t \right)}^{\chi (X)}}
$$
The present article,
in particular, recovers the result for
$X\,=\,{\bf P}^{n}$ which was worked out in \cite{la&yau-hosy}.

 \section{Varieties with a Torus Action}

Throughout this section $X$ means a
projective algebraic variety, on which an algebraic torus $T$  acts linearly
 having only a finite number of
invariant irreducible subvarieties of dimension p.
 In particular, we will see that
the result is true for any projective toric variety. Let us denote by $H$
 the homology group $H_{2p}\,(X,{\bf Z}).$

The action of $T$ on $X$ induces an action on the Chow variety
 ${\cal C} _{p,d} (X).$
Let $\lambda$ be an element in $H$ and denote by ${\cal C}_{\lambda}^{T}$
the fixed point set of ${\cal C} _{\lambda} (X)$ under the action of $T.$
Then its Euler characteristic
$\chi \left( {\cal C}_{\lambda}^{T} \right)$ is equal to the number
of invariant subvarieties of $X$ with homology class $\lambda.$
We have
\begin{equation}
\label{lawson}
E_{p}(\lambda) \, = \, \chi \left( {\cal C}_{\lambda} \right) \, =
                      \, \chi \left( {\cal C}_{\lambda}^{T} \right)
\end{equation}
where the  first equality is just the definition of $E_{p}$
and the last one is proved in  \cite{la&yau-hosy}. The following
theorem tells us that  $E_{p}$ is rational.
\begin{theorem}
\label{main}
Let $E_{p}$ be the Euler  series of $X.$  Denote by
$V_1, \ldots , V_N$ the p-dimensional invariant irreducible
subvarieties of $X.$ Let $e_{[V_i]} \in {\bf Z}[C]$
be the characteristic function of the subset $\{ [V_i] \}$ of $C.$
Then
\begin{equation}
\label{ecp}
E_{p}=\, \prod_{1\leq i\leq N}\left( \frac{1}{1-e_{[V_i]}}\right)
\end{equation}
\end{theorem}
{\bf Proof:}
 For each $V_i$ define $f_i$ in ${\bf Z}[[C]]$ by
\begin{equation}
\label{function}
f_i(\lambda)=\left\{ \begin{array}{ll} 1 & \mbox{if $\lambda=n\cdot
                                               [V_i]$, $n\geq 0$} \\
                                      0 & \mbox{otherwise.}\end{array} \right.
\end{equation}
It is easy to see from equations (\ref{lawson}) and (\ref{function})
 that $E_{p}$ can be written as
\begin{equation}
\label{ecr}
E_{p} \, = \, \prod_{1\leq i\leq N} f_i.
\end{equation}
 The theorem follows because of the equality
$$
1=\left( 1-e_{[{V}_{i}]} \right) \cdot f_i
$$

\qed

Observe that if we fix a basis for  $H$ modulo torsion, then
the elements of ${\bf Z}[C]$ can be identified  with some
Laurent polynomials. Under this identification
 any rational element  of
${\bf Z}[[C]]$ is a rational function.

The next lemma tells us that the result is true for
any projective toric variety.

\begin{lemma}
\label{ciclos}
Let $X$ be a projective (perhaps singular) toric variety. Then any
irreducible subvariety V of $ X$  which is
invariant under the torus action is the closure
of an orbit. Therefore, any invariant cycle has the form
$$
c \, = \, \sum \, n_{i} \, \overline{\cal O}_{i}
$$
where
each $n_{i}$ is a nonnegative integer and each $\overline{\cal O}_{i}$
is the closure of the orbit ${\cal O}_{i}$.

\end{lemma}

{\bf Proof: }
 The fan $\Delta$ associated to $X$ is finite because $X$
 is compact. Hence there is a finite number
of cones, therefore a finite number of orbits.
Let V be an invariant irreducible subvariety of $X.$ We can express V as
 the closure
of the union of orbits. Since there is a finite number of them
 we must have that
$$
V \, = \, {\overline{\cal O}}_{1} \, \cup \, {\overline{\cal O}}_{2}
 \, \cup \cdots  \cup \, {\overline{\cal O}}_{N}
$$
where $\overline{\cal O}_{i}$ is the closure of an orbit. Finally
 since V is irreducible, there must be
$i_{0}$ such that  $V \, = \, {\overline{\cal O}}_{{i}_{0}}$.
\qed

\section{Smooth Toric Varieties}

In this section we give an equivariant version of the
Euler series and find a relation between the equivariant
and not equivariant Euler series.
The use of equivariant cohomology allows us to analize the Euler
series from
a geometrical point of view. This approach
might help to understand other cases.
Throughout this section, unless otherwise stated,
 $X$ is a smooth projective toric variety, and
we use cohomology instead of
homology by applying Poincar\'e duality.

Let $H$ and $H_{T}$ be the cohomology group $H^{2p}(X,{\bf Z})$
and the equivariant cohomology group $H_{T}^{2p}(X,{\bf Z})$ of $X,$
respectively. Denote by $\Delta$ the fan associated with $X.$

Let ${\cal C}_{\lambda}^{T}$ and ${\cal C}_{\lambda}$ be the spaces of
 all p-dimensional effective invariant cycles
and p-dimensional effective cycles on $X$ with cohomology class $\lambda.$
 It is proved in \cite{la&yau-hosy}
that
\begin{equation}
\label{blaine}
\chi  \left( {\cal C}_{\lambda}^{T} \right) \, = \, \chi
 \left( {\cal C}_{\lambda} \right).
\end{equation}
The next lemma is crucial for the following results,
\begin{lemma}
\label{finito}
Let $\lambda$ be an element in $H.$ Then
${\cal C}_{\lambda}^{T}$ is a finite set.
\end{lemma}
{\bf Proof:}
 By lemma ~\ref{ciclos} we know that any
invariant effective cycle $c$  in  ${\cal C}_{\lambda}^{T}$
has the form
 $c \, = \, \sum_{i=1}^{N} \, {\beta}_{i} \, {\overline{\cal O}}_{i} $
 with $\beta_{i} \in {\bf N},$ we obtain  that
 ${\cal C}_{\lambda}^{T} $
 has a countable number of elements. We know that
 ${\cal C}_{\lambda}$ is a projective algebraic
variety (see lemma ~\ref{compacto}), and since
${\cal C}_{\lambda}^{T}$
 is Zariski closed in ${\cal C}_{\lambda},$ we have that
${\cal C}_{\lambda}^{T}$ is a finite set.
\qed

Our next step is to define the equivariant Euler series
for $X.$ \\
Let
${\overline{\cal O}}$ be an irreducible invariant cycle in a smooth toric
 variety (lemma ~\ref{ciclos}). Since
${\overline{\cal O}} \subset X$ is smooth, we have an equivariant
 \mbox{Thom-Gysin}
sequence
$$
\cdots \longrightarrow
H_{T}^{i-2cod\,\overline{\cal O}} (\overline{\cal O})
\longrightarrow \, H_{T}^{i} \, (X)
\longrightarrow \, H_{T}^{i} \,(X-\overline{\cal O})
\longrightarrow \, \cdots
$$
and we define ${[\overline{\cal O}]}_{T}$ as the image of 1 under
$$
H_{T}^{0} \, (\overline{\cal O})
\longrightarrow H_{T}^{2cod\,\overline{\cal O}} \, (X).
$$
Let $\{ D_{1}, \ldots, D_{K} \}$ be the set of T-invariant divisors on $X.$
To each $D_{i}$ we associate the variable \, $t_{i}$ \,
in the polynomial ring
\,${\bf Z}[t_{1}, \ldots, t_{K}].$\,
Let $\cal I$ be the ideal generated by the (square free) monomials
$\{ t_{{i}_{1}} \cdots  t_{{i}_{l}} | \,
 D_{{i}_{1}}  + \cdots  + D_{{i}_{l}} \not\in \, \Delta \}.$
It is proved in \cite{bcp-regular} that
\begin{equation}
 \label{eqring}
{\bf Z}[t_{1}, \ldots, t_{K}] / {\cal I} \, \cong \,
H_{T}^{\ast} (X, \, {\bf Z})
\end{equation}
The arguments given there also prove the following.
\begin{proposition}
\label{independent}
For any T-orbit $\cal O$ in a smooth projective toric variety $X$, one has
$$
[\overline{\cal O}]_{T} \, = \, \prod_{{\overline{\cal O}}
\subset{D_{i}}} \, [D_{i}]_{T} .
$$
Furthermore if $\cal O$ and ${\cal O}^{\prime}$ are distinct orbits, then
$$
[{\overline{\cal O}}]_{T} \, \not= \,  [\overline{{\cal O}^{\prime}}]_{T}.
$$
\end{proposition}
It is natural to define
 the cohomology class for any effective invariant cycle \
$V = \sum \, m_{i} \, {\overline{\cal O}}_{i}$ \  as \
${[V]}_{T} = \sum_{i} \, m_{i} \, {[{\overline{\cal O}}_{i}]}_{T}.$
 where ${\cal O}_{i} \not= \, {\cal O}_{j}$\ if $i\not= j$
In a similar form as we define $C,$ ${\bf Z}[[C]]$ and
${\bf Z}[C],$
we denote by $C_{T}$ the monoid of equivariant cohomology
classes of invariant effective cycles of codimension p,
by ${\bf Z}[[C_{T}]]$  the set of functions on $C_{T},$
and by ${\bf Z}[C_{T}]$ the set of functions with finite support on $C_{T}.$
Since $C_{T} \simeq {\bf N}^{N}$ where $N$ is the number of orbits of
codimension p, we obtain that
$$
+: \, C_{T} \times C_{T} \, \longrightarrow \, C_{T}
$$
has finite fibers.
Observe that if $\pi : H_{T} \rightarrow H$ denotes the standard
surjection, we obtain from lemma ~\ref{finito} that
$$
\pi : C_{T} \longrightarrow C
$$
 is onto with finite fibers.
We arrive at the following definition:
\begin{definition}
\label{eec}
Let $X$ be a smooth projective toric variety and let
$H_{T} \, (X)$ be the equivariant cohomology of $X$.
 Let us denote by ${\cal C}^{T}_{\xi}$
the space of all invariant effective cycles on $X$ whose
 equivariant cohomology class is $\xi .$
The {\bf equivariant Euler series} of $X$ is the element
$$
E_{p}^{T} = \sum_{\xi} \, \chi \left( {\cal C}^{T}_{\xi} \right)\,  \xi
 \; \; \; \in \; \; {\bf Z}[[C_{T}]]
$$
 where  the sum is over $C_{T}.$
\end{definition}
Let us define the  ring homomorphism
$$
J: \, {\bf Z}[[C_{T}]] \, \longrightarrow \,
             {\bf Z}[[C]]
$$
by
$$
J(\xi) \, = \,
      \sum_{\lambda} \left( \sum_{\pi(\beta)=\lambda} \, a_{\beta}
             \right) \lambda
$$
where $\xi \, = \, \sum_{\beta} \, a_{\beta}\, \beta.$ This is
well defined since $\pi$ has finite fibers.

\begin{theorem}
\label{ec-eec}
Let $X$ be a smooth projective toric variety. Denote by
$E_{p}$, $E_{p}^{T}$ and $J$ the Euler series,
the equivariant Euler series and the ring homomorphism
defined above. Then
$ J \left( E_{p}^{T} \right) \, = \, E_{p}.$
Furthermore,
$$
E_{p}^{T} \,=\,\prod_{1\leq{i}\leq{N}} \,
   \left( \frac{1}{1-e_{[{\overline{\cal O}}_i]_{T}}} \right)
$$
and therefore
$$
E_{p} \,=\,\prod_{1\leq{i}\leq{N}} \,
   \left( \frac{1}{1-e_{[{\overline{\cal O}}_i]}} \right)
$$
\end{theorem}
{\bf Proof: }
We define for each orbit ${\cal O}_{i}$ an element
$ f_{i}^{T} \in {\bf Z}[[C_{T}]]$ by
\begin{equation}
\label{ec.e}
f^T_i(\xi)=\left\{ \begin{array}{ll}
1 & \mbox{if $\xi=n\cdot
                  [{\overline {\cal O}_i}]_T$, $n\geq 0$} \\
0 & \mbox{otherwise.\mbox{              } } \end{array} \right.
\end{equation}
and denote by $e_{\xi}$ the characteristic function of $\{\xi\}.$
It follows from  both definition
 (\ref{eec}) and equation (\ref{ec.e}) that
$$
E_{p}^{T} \, = \, \prod_{1\leq{i}\leq{N}} \, f_{i}^{T} \,
$$
and
$$
\left( 1-e_{[{\overline {\cal O}_{i}}]_T} \right)
                                \cdot f_{i}^{T} \, = \, 1.
$$
Therefore the equivariant Euler series is rational and
\begin{equation}
\label{rationale}
E_{p}^{T} \, = \,\prod_{1\leq{i}\leq{N}}
    \left(  \frac{1}{1-e_{[{\overline{\cal O}_{i}}]_{T}}} \right).
\end{equation}
For each $V_i$ we defined (see equation \ref{function})
a function $f_i$ in ${\bf Z}[[C]]$ by
$$
f_i(\lambda)=\left\{ \begin{array}{ll} 1 & \mbox{if $\lambda=n\cdot
                                               [V_i]$, $n\geq 0$} \\
                                      0 & \mbox{otherwise.}\end{array} \right.
$$
And we know from  theorem ~\ref{main} that
$$
E_p \, = \, \prod_{i=1}^{N} \, f_i
\mbox{  \, \, \, \, \, \,  with \, \, \, \, \, \,}
f_i \cdot \left( 1-e_{[V_{i}]} \right) \, = \, 1.
$$
Now, the result follows since
$\pi ([{\overline{\cal O}}_i]_{T}) \, = \,
[{\overline{\cal O}}_{i}]$ and  $J$ is a ring homomorphism satisfying
$$
J \left( f_{i}^{T} \right) \,=\, f_{i}
\mbox{ \, \, \, \, \, \,    and  \, \, \, \, \, \,  }
J \left( e_{\sigma}\right) \,=\, e_{\pi(\sigma)}.
$$
\qed

\section{Some Examples}
\noindent {\bf  I) The projective space ${\bf P}^{n}$}

Let  \, $X = {\bf P}^{n}$ \, be the complex projective space of dimension n.
Let $\{e_{1},\ldots,e_{n}\}$ be the standard
basis for ${\bf R}^{n}$. Consider $A \, = \,
 \{e_{1},\ldots,e_{n+1}\}$ a set of generators of the
 fan $\Delta$ \,
where $e_{n+1} = - \sum_{i=1}^{n} e_{i}$. We have the following equality
 $$
H^{\ast}(X, \, {\bf Z}) \, \cong  \, {\bf Z}
 \, [t_{1}, \ldots, t_{n+1}] \, /I
$$
where I is the ideal generated by
 $$
i) \; \; \;   t_{1} \cdots t_{n+1}
$$
and
$$
ii) \; \; \;   \sum_{j=1}^{n+1} \, e_{i}^{\ast} (e_{j}) \,  t_{j} \; \;\;
 \mbox{\, \, for\, \,  }i\,=\,1,\ldots,n \, ,
$$
where $e_{i}^{\ast} \, \in \, ({\bf R}^{n})^{\ast}$
is the element dual to $e_{i}$.\\
However $ii)$ says that   $\; t_{i} \sim t_{j} \; $ for all. Therefore
$$
H^{\ast} \, (X, \, {\bf Z}) \, = \, {\bf  Z} \, [t] \, / t^{n+1}.
$$
Consequently, any two cones of dimension p represent the same element in
 cohomology, and
Theorem  ~\ref{ec-eec} implies that
$$
{\displaystyle \prod_{i=1}^{(_{\; \; \, p}^{n+1})} {\left( \frac{1}{1-t}
\right)}
 \, = \,
\left( \frac{1}{1-t} \right)^{(_{\; \; \; p}^{n+1})} \, = \, E_{p}}.
$$

\noindent {\bf II) $P^{n} \times P^{m}$}

Recall that $X( \Delta \times {\Delta}^{\prime}) \,
                  \cong \, X(\Delta) \times X({\Delta}^{\prime})$.
 Using the same
notation as in example I, we have that a set of generators of
 $\Delta \times \Delta '$  is given by
$$
\{e_{1}, \, \ldots, \, e_{n}, e_{n+1}, \ldots, e_{n+m}, e_{n+m+1}, e_{n+m+2}\}
$$
with $e_{n+m+1} \, = \, - \sum_{i=1}^{n} \, e_{i}$ and $e_{n+m+2} \, = \,
 - \sum_{i=n+1}^{n+m} \, e_{i}.$ and
 $\{e_{1}, \ldots, e_{n+m}\}$ is a basis for $P^{n} \times P^{m}$. Then
$$
H^{\ast} \, (X, \, {\bf Z}) \, = \, {\bf Z} \, [t_{1}, \ldots , t_{n+m+2}] \, /
\, I
$$
where I is  the ideal generated by
$$
i)\ \  \{t_{1} \cdots t_{n} t_{n+m+1}, \; \; t_{n+1}
\cdots t_{n+m}t_{n+m+2}, \; \;
 \prod_{i=1}^{n+m+2} t_{i} \}
$$
and
$$
ii)\  \ \sum_{j=1}^{n+m+2} \, e_{i}^{\ast} (e_{j}) \,  t_{j}\  \; \; \; \; i\,
 = \, 1,\ldots, n+m \, .
$$
{}From $ii)$  we obtain,
\begin{eqnarray*}
t_{i} \, \sim \, t_{n+m+1} & & \mbox{if \ \  $1 \leq i \leq n $} \\
t_{j} \, \sim \, t_{n+m+2} & & \mbox{if \ \  $n+1 \leq j \leq n+m $}
\end{eqnarray*}
The number of cones of dimension p is equal to $\sum_{k+l=p} \, (_{\; \; \;
k}^{n+1})
 \,(_{\; \; \; l}^{m+1})$. \\
Denote by ${\displaystyle t_{k,l} \, = \, t_{{i}_{1}} \cdots t_{{i}_{k}}
t_{{j}_{1}}
 \cdots t_{{j}_{l}}}.$ Then
$$
\prod_{k+l=p} {\left( \frac{1}{1-t_{k,l}} \right)}^{(_{\; \; \,k}^{n+1}) \,
 (_{\;\;\,l}^{m+1})} \, = \, E_{p}
$$
where $t_{k,l} \, = \, t_{n+m+1}^{k} \, t_{n+m+2}^{l}$.

\noindent {\bf III) Blow up of ${\bf P}^{n}$ at a point.}

The fan $\tilde{\Delta}$ associated to the blow up $\tilde{{\bf P}^{n}}$
 of the projective space
at the fixed point given by the cone ${\bf R}^{+} e_{2} +
 \cdots + {\bf R}^{+} e_{n+1}$ is
generated by $\{e_{1}, \, \ldots, \, e_{n+1}, e_{n+2}\}$ where $e_{n+2} \,
 = \, -e_{1}$.
Denote by $D_{i}$ the 1-dimensional cone \, ${\bf R}^{+}\, e_{i}$ \,
 and by $s_{i}$  its class in cohomology where
$$
H^{\ast} \, (X, \, {\bf Z}) \, = \, {\bf Z} \, [s_{1}, \ldots , s_{n+2}] \, /
\, I
$$
and I is the ideal generated by
$$
i)\ \  \{s_{{i}_{1}} \cdots s_{{i}_{k}}\, |  \, D_{{i}_{1}}+
\cdots + D_{{i}_{k}}
 \, \mbox{is not in } \tilde{\Delta} \}
$$
and
$$
ii)\  \ \sum_{j=1}^{n+2} \, e_{i}^{\ast} (e_{j}) \,
 s_{j}\  \; \; \; \; i\, = \, 1,\ldots, n \, .
$$
However $ii)$ is equivalent to
$$ii) \, \, \, s_{2}\sim \cdots \sim s_{3} \sim s_{n+1}
              \; \; \; \mbox{and} \; \; \;
 s_{1} \sim s_{n+1}+s_{n+2}$$
Note that a p-dimensional cone cannot contain
both $D_{n+2}$ and  ${D_{1}}$. The
 reason is  that $D_{n+2}$ is
generated by $- e_{1}$ and  $D_{1}$ by $e_{1},$ but by definition,
 a cone does not contain a subspace of
dimension greater than 0.
 We would like to find a basis for $H^{\ast}(\tilde{{\bf P}^{n}})$ and
 write any monomial
of degree p in terms of it.
Consider the monomial $ s_{{i}_{1}} \cdots  s_{{i}_{p}}$.
 There are three possible situations:

\noindent {\bf 1)}  $s_{{i}_{j}}$ is different from both  $s_{n+2}$ and
$s_{n+1}$. In
this situation we have from $ii)$ that
$s_{{i}_{1}} \cdots  s_{{i}_{p}} \, = \, s_{n+1}^{p}$.

\noindent{\bf 2)}  $s_{n+2}$ is equal to $s_{{i}_{j}}$ for some $j \, =
 \, 1,\ldots, p$. Then from $ii)$
we obtain that $s_{{i}_{1}} \cdots  s_{{i}_{p}} \, = \, s_{n+1}^{p-1} \,
s_{n+2}$.

\noindent {\bf 3)}  $s_{1}$ is equal to $s_{{i}_{j}}$ for some $j \, = \,
1,\ldots, p$.
 Then from $ii)$ we obtain
$s_{{i}_{1}} \cdots  s_{{i}_{p}} \, = \, (s_{n+1} + s_{n+2}) \, s_{n+1}^{p+1}
\, = \, s_{n+1}^{p} + s_{n+2}  \, s_{n+1}^{p-1}$ which is the sum of 1) and 2).

\noindent We conclude that $s_{n+1}^{p}$ and $s_{n+2} s_{n+1}^{p-1}$ form a
 basis for $H^{2p}$ if $p < n$.
If $p=n$ then $s_{n+1}^{p} \, = \, 0$ and the only generator is
 $s_{n+2} s_{n+1}^{p-1}$. Let us call
$s_{n+1}$ by $t_{1}$ and $s_{n+2} s_{n+1}^{p-1}$ by $t_{2}$.
 The Euler series for
 $\tilde{{\bf P}^{n}}$ is:
$$
{\displaystyle E_{p} \, = \, {\left(\frac{1}
{1-t_{1}}\right)}^{(_{p}^{n})} \,
 {\left(\frac{1}{1-t_{1}t_{2}}\right)}^{(_{p-1}^{\, \, \, n})}
\, {\left(\frac{1}{1-t_{2}}\right)}^{(_{p-1}^{\, \, \, n})}}
 \; \; \; \;{\mbox if }\; \; p < n \, .
$$
$$
{\displaystyle E_{p} \, = \,
 {\left(\frac{1}{1-t_{2}}\right)}^{(_{\; \; \; p}^{n+2})}}
 \; \; \; \;{\mbox if }\; \; p\,=\,n \, .
$$

\noindent {\bf IV) Hirzebruch surfaces}

A set of generators for the fan $\Delta$
 that represents the Hirzebruch surface $X(\Delta)$ is
given by $\{e_{1}, \ldots , e_{4} \}$ with
$\{e_{1}, e_{2}\}$ the standard basis for ${\bf R}^{2}$, and
 $e_{3} \, = \, -e_{1} + ae_{2}, \; \; a > 1$
and $e_{4} \, = \, -e_{2}$.
 With the same notation as in the last examples, we have
$$
H^{\ast} (X(\Delta)) \, = \, {\bf Z}[t_{1}, \ldots, t_{4}] \, / \, I
$$
where I is generated by
$$
 i) \; \; \; \{ t_{1}t_{3}, \, t_{2}t_{4} \}
$$
and
$$
ii) \; \; \; \{ t_{1} -t_{3}, \, t_{2}+ at_{3} -t_{4} \}
$$
from $ii)$ we have the following conditions for the $t_{i}$'s in $H^{\ast}(X)$
\begin{equation}
\label{coh}
  \; \; \; \; t_{1} \, \sim \, t_{3} \; \; \mbox{and }\; \; t_{2}
 \, \sim \, (t_{4}-at_{3}).
\end{equation}
A basis for $H^{\ast} (X)$ is given by $\{ \{0\}, t_{3}, t_{4}, t_{4}t_{1}  \}$
 \, (see \cite{da-tova},\cite{ful-tova}).
The Euler series for each dimension is:

\noindent  {\bf 1)} Codimension 2:
 There are four orbits (four cones of dimension 2), and all
of them are equivalent in homology.  From Theorem  ~\ref{ec-eec} we obtain
$$
E_{2} \, = \, {\left(\frac{1}{1-t}\right)}^{4}
$$
 {\bf 2)}  Codimension 1:
 Again, there are four orbits (four cones of dimension 1), and the
 relation  among them,
in homology, is given by ~\ref{coh}. From theorem ~\ref{ec-eec}  we obtain
$$
E_{1} \, = \, {\left(\frac{1}{1-t_{3}}\right)}^{2} \,
 \left(\frac{1}{1-t_{4}}\right) \, \left(\frac{1}{1-t_{3}^{-a} t_{4}}\right).
$$
 {\bf 3)}  Codimension 0: The only orbit is the torus itself so
$$
 E_{0} \, = \, \left(\frac{1}{1-t}\right).
$$

\bibliography{referenc}

\begin{thebibliography}{BDP90}

\bibitem[BDP90]{bcp-regular}
E.~Bifet, C.~{De Concini}, and C.~Procesi.
\newblock Cohomology of {R}egular {E}mbeddings.
\newblock {\em Adv. in Math.}, 82:1--34, 1990.

\bibitem[Bif92]{emili-cartas}
E.~Bifet.
\newblock Personal letters.
\newblock August 1 and December 5, 1992.

\bibitem[Dan78]{da-tova}
V.~I. Danilov.
\newblock The theory of toric varieties.
\newblock {\em Russian Math. Surveys}, 33(2):97--154, 1978.

\bibitem[Eli92]{eli-thesis}
J.~Elizondo.
\newblock {\em The Euler-Chow Series for Toric Varieties}.
\newblock Ph.D. thesis, SUNY at Stony Brook, August 1992.

\bibitem[Ful93]{ful-tova}
W.~Fulton.
\newblock {\em Introduction to Toric Varieties}.
\newblock Number 131 in Annals of Mathematics Studies. Princeton University
  Press, Princeton, NJ, 1st. edition, 1993.

\bibitem[LY87]{la&yau-hosy}
H.~B. Lawson, Jr. and S.~S.~T. Yau.
\newblock Holomorphic symmetries.
\newblock {\em Ann. scient. \'{E}c. Norm. Sup.}, t. 20:557--577, 1987.


\bibitem[Mac62]{mcd-sym}
I.~G. MacDonald.
\newblock The {P}oincar\'{e} polynomial of a symmetric product.
\newblock {\em Proc. Cam. Phil. Soc.}, 58:563--568, 1962.

\end{thebibliography}

\bibliographystyle{alpha}

\begin{minipage}[t]{2in}
{\footnotesize Instituto de Matem\'aticas \vspace{-0.1cm}\newline
Ciudad Universitaria, UNAM \vspace{-0.1cm}\newline
M\'exico D.F. 04510} \vspace{-0.1cm}\newline
M\'exico
\end{minipage}
\ {\footnotesize and} \
\begin{minipage}[t]{2in}{\footnotesize
Mathematics Section, ICTP.\vspace{-0.1cm}\newline
P.O.Box 586  \vspace{-0.1cm}\newline
34100 Trieste \vspace{-0.1cm}\newline
Italy \newline\vspace{-0.2cm}}
\end{minipage}
\newline
{\footnotesize  e-mail: elizondo@ictp.trieste.it}

\end{document}